\documentclass[12pt]{article}
\usepackage{amsmath, amssymb}
\usepackage{amsfonts}
\usepackage{graphicx}
\usepackage{latexsym}
\usepackage{url}
\usepackage{a4wide,amsfonts,amsmath,latexsym,amssymb,euscript,graphicx,units,mathrsfs}
\usepackage{graphicx}
\usepackage{color}
\usepackage{amssymb}
\usepackage{amssymb}
\usepackage[T1]{fontenc}
\usepackage{latexsym}
\usepackage{xypic}
\usepackage{eufrak}
\usepackage{euscript}
\usepackage{amsfonts,amsmath}
\usepackage{verbatim}
\usepackage{fancyhdr}
\usepackage[english]{babel}
\usepackage{mathrsfs}
\usepackage{units}

\newcommand{\qed}{$\hfill\Box$}

\newcommand{\R}{\mathbb{R}}

\newcommand{\E}{\mathbb{E}}
\newcommand{\PP}{\mathbb{P}}
\newcommand{\N}{\mathbb{N}}
\newcommand{\Z}{\mathbb{Z}}

\newcommand{\HH}{\mathfrak{H}}

\def\sk{{\mathbb{D}}}
\def\rcr{\right]}
\def\lcr{ \left[ }

\newcommand{\schw}{\stackrel{d}{\longrightarrow}}

\newcommand{\bee}{\begin{equation}}
\newcommand{\eee}{\end{equation}}
\newcommand{\bea}{\begin{eqnarray}}
\newcommand{\eea}{\end{eqnarray}}
\newcommand{\bean}{\begin{eqnarray*}}
\newcommand{\eean}{\end{eqnarray*}}

\renewcommand{\theequation}{\arabic{section}.\arabic{equation}}

\newtheorem{prop}{Proposition}[section]
\newtheorem{cor}[prop]{Corollary}

\newtheorem{lem}[prop]{Lemma}
\newtheorem{ex}[prop]{Example}
\newtheorem{theo}[prop]{Theorem}

\begin{document}

\title{\bf Quantitative Breuer-Major Theorems}
\author{Ivan Nourdin\thanks{Laboratoire de Probabilit\'es et Mod\`eles Al\'eatoires, Universit\'e Pierre et Marie Curie, Bo\^ite courrier 188, 4
Place Jussieu, 75252 Paris Cedex 5, France, Email: ivan.nourdin@upmc.fr.} \and
Giovanni  Peccati\thanks{Facult\'{e} des Sciences, de la Technologie
et de la Communication; UR en Math\'{e}matiques. 6, rue Richard
Coudenhove-Kalergi, L-1359 Luxembourg, Email: giovanni.peccati@gmail.com.} \and
Mark Podolskij\thanks{Department of Mathematics, ETH Z\"urich, HG G32.2, 8092 Z\"urich,
Switzerland, Email: mark.podolskij@math.ethz.ch.}}

\date{June 3, 2010}

\maketitle

\begin{abstract}
We consider sequences of random variables of the type $S_n= n^{-1/2} \sum_{k=1}^n 
\{f(X_k)-E[f(X_k)]\}
$, $n\geq 1$, where $X=(X_k)_{k\in \Z}$ is a $d$-dimensional
Gaussian process and $f: \R^d \rightarrow \R$ is a measurable function.
It is known that, under certain conditions on $f$ and the covariance function $r$ of $X$, 
$S_n$ converges in distribution to a normal variable $S$. In the present paper we derive several explicit 
upper bounds for quantities of the type $|\E[h(S_n)] -\E[h(S)]|$, where $h$ is a sufficiently smooth test 
function. Our methods are based on Malliavin calculus, on interpolation techniques and on the Stein's method 
for normal approximation. The bounds deduced in our paper depend only on 
${\rm Var}[f(X_1)]$ 
and on simple infinite 
series involving the components of $r$. In particular, our results generalize and refine some classic CLTs by 
Breuer-Major, Giraitis-Surgailis and Arcones, 
concerning the normal approximation of partial sums associated with Gaussian-subordinated time-series.

\end{abstract}




\eject

\section{Introduction} \label{Intro}
\setcounter{equation}{0}
\renewcommand{\theequation}{\thesection.\arabic{equation}}

Fix $d\geq 1$ and consider a $d$-dimensional centered stationary Gaussian process $X=(X_k)_{k\in \Z}$, $X_k=(X_k^{(1)}, \ldots, X_k^{(d)})$,
defined on a probability space $(\Omega, \mathcal F, \PP)$. For any $1\leq i,l\leq d$ and $j\in \Z$, we denote by
\begin{equation} \label{corr}
r^{(i,l)}(j) = \E[X^{(i)}_1 X^{(l)}_{1+j}],
\end{equation}
the covariance of $X^{(i)}_1$ and $X^{(l)}_{1+j}$. Let $f: \R^d \rightarrow \R$ be a measurable function
and write
\begin{equation} \label{sn}
S_n= \frac{1}{\sqrt n} \sum_{k=1}^{n} 
\{f(X_k)-\E[f(X_k)]\}, 
\quad n\geq 1,
\end{equation}
to indicate the sequence of {\sl normalized partial sums} associated with the subordinated process $k\mapsto f(X_k)$. One crucial problem in Gaussian analysis is the following:

\begin{quotation}
\noindent {\bf Problem P.} {\sl Find conditions on $f$ and on the covariance $r$ in order to have that, as $n\rightarrow\infty$, $S_n$ converges in distribution to a Gaussian random variable. }
\end{quotation}

Albeit easily stated, Problem P is indeed quite subtle. For instance, as observed e.g. in \cite[p. 429]{BM}, it is in general not possible to deduce a solution to Problem P by using standard central limit results for dependent random variables (for instance, by applying techniques based on mixing). More to the point, slight variations in the form of $f$ and $r$ may imply that either the normalization by the factor $n^{1/2}$ is inappropriate, or the limiting distribution is not-Gaussian (or both): see Dobrushin and Major \cite{DobMaj}, Rosenblatt \cite{rosnblatt1961} and Taqqu \cite{Taqqu1, Taqqu2} for several classic results connected to this phenomenon, as well as Breton and Nourdin \cite{breton-nourdin} for recent developments.

It turns out that an elegant solution to Problem P can be deduced by using the notion of {\sl Hermite rank}.
Recall that the function $f$ is said to have Hermite rank equal to $q$ with respect to $X$, 
where $q\geq 1$ is an integer, if { (a)} $\E[(f(X)-\E[f(X)])p_m(X)]=0$ for every polynomial $p_m$ (on $\R^d$) 
of degree $m\leq q-1$; and { (b)} there exists a polynomial $p_q$ of degree $q$ such that 
$\E[(f(X)-E[f(X)])p_q(X)]\neq 0$ 
(see also Proposition \ref{D : HermiteRank} below). Then, one has the following well-known statement:

\begin{theo}[Breuer-Major theorem for stationary vectors] \label{th1}
Let $\E[f^2(X_1)]<\infty$, and assume that the function $f$ has Hermite rank equal to $q\geq 1$. Suppose that
\begin{equation}\label{zoe}
\sum_{j\in \Z} |r^{(i,l)}(j)|^q<\infty, \quad \forall  i,l\in\{1,...,d\}.
\end{equation}
Then 
$\sigma^2:= {\rm Var}[f^2(X_1)] +2\sum_{k=1}^\infty {\rm Cov}[f(X_1),f(X_{1+k})]$ 
is well-defined, and belongs to $[0,\infty)$. Moreover, one has that
\begin{equation} \label{clt}
S_n \schw S\sim N(0,\sigma^2),
\end{equation}
where $N(0,\sigma^2)$ indicates a centered Gaussian distribution with variance $\sigma^2$, and $\schw$ stands for convergence in distribution.
\end{theo}
In the case $d=1$, Theorem \ref{th1} was first proved by Breuer and Major in \cite{BM}, whereas Theorem 4 in Arcones \cite{AR} proves the statement for a general $d$ (both proofs in \cite{AR, BM} are based on the method of cumulants and diagram formulae -- see e.g. \cite{PecTaqBOOK, Sur}).  The reader is referred to Sun \cite{SUN1965} for an early statement in the case of a Hermite rank equal to 2, and to Giraitis and Surgailis \cite{GS} for some continuous-time analogous of Theorem \ref{th1}. Note that any central limit result involving Hermite ranks and series of covariance coefficients is customarily called a `Breuer-Major Theorem', in honor of the seminal paper \cite{BM}.

Theorem \ref{th1} and its variations have served as fundamental tools for Gaussian approximations in an impressive number of applications, of which we provide a representative (recent) sample: renormalization of fractional diffusions \cite{ahnLeon}, power variations of Gaussian and Gaussian-related continuous-time processes \cite{BCP, BCP2, CoNuWo, LeonLudena2007}, Gaussian fluctuations of heat-type equations \cite{BeghnAL}, estimation of Hurst parameters of fractional processes \cite{BrezLeon, coeurjolly2005, CoeurAoA2008}, unit root problems in econometrics \cite{BuChan}, empirical processes of long-memory time-series \cite{Mar2002}, level functionals of stationary Gaussian fields \cite{KtatzLeo}, variations of multifractal random walks \cite{LudAAP2008}, and stochastic programming \cite{Wang}. See also Surgailis \cite{Sur} for a survey of some earlier uses of Breuer-Major criteria.

Despite this variety of applications, until recently the only available techniques for proving results such as Theorem \ref{th1} were those based on combinatorial cumulants/diagrams computations. These techniques are quite effective and flexible (see e.g. \cite{MarPTRF, sodinTS} for further instances of their applicability), but suffer of a fundamental drawback, namely they {\sl do not} allow to deduce Berry-Esseen relations of the type
\begin{equation}\label{TEST1}
\big|\E[h(S_n)] - \E[h(S)]\big|\leq \varphi(n), \quad n\geq 1,
\end{equation}
where $h$ is a suitable test function, and $\varphi(n)\rightarrow 0$ as $n\rightarrow \infty$. An upper bound such as (\ref{TEST1}) quantifies the error one makes when replacing $S_n$ with $S$ for a fixed $n$.

In \cite[Section 4]{NP1}, the first two authors of this paper proved that one can combine Malliavin calculus (see e.g. \cite{N}) and Stein's method (see e.g. \cite{Chen_Shao_sur}) to obtain relations such as (\ref{TEST1}) (for some explicit $\varphi(n)$) in the case where: (i) $d=1$, (ii) $f=H_q$ is a Hermite polynomial of degree $q\geq 2$ (and thus has Hermite rank equal to $q$), (iii) $X$ is obtained from the increments of a fractional Brownian motion of Hurst index $H<1-(2q)^{-1}$, and (iv) $h$ is either an indicator of a Borel set or a Lipschitz function. Since under (iii) one has that $|r(j)|\sim j^{2H-2}$, these findings allow to recover a very special case of Theorem \ref{th1} (see Example
\ref{fBM} below for more details on this point).

The aim of the present work is to extend the techniques initiated in \cite{NP1} in order to deduce several complete {\sl quantitative Breuer-Major theorems}, that is, statements providing explicit upper bounds such as (\ref{TEST1}) for any choice of $f$ and $r$ satisfying the assumptions of Theorem \ref{th1}. We stress by now that we will {\sl not} require that the functions $f$ enjoy any additional smoothness property, so that our results represent a genuine extension of the findings by Breuer-Major and Arcones.

As anticipated, our techniques are based on the use of Malliavin operators on a Gaussian space, that we combine both with Stein's method and with an interpolation technique (already applied in \cite{NouPeReinertAOP, PecZheng}) which is reminiscent of the `smart path method' in Spin Glasses -- see e.g. Talagrand \cite{talag}. In particular, the use of Stein's method allows to deal with functions $h$ that are either Lipschitz or indicators of the type $h = {\bf 1}_{(-\infty,z]}$, whereas the use of interpolations requires test functions that are twice differentiable and with bounded derivatives. Note that this implies that the convergence (\ref{clt}) takes indeed place in the stronger topologies of the Kolmogorov and Wasserstein distances.

The rest of the paper is organized as follows. Section \ref{Prelim} contains the statements of our main results, some examples
and applications. Section \ref{S : toolbox} presents some notions and results that are needed to prove our main Theorem \ref{mainth}. Section \ref{S : PROOFS} is devoted to proofs.

\section{Statement of the main results}\label{Prelim}
We keep the assumptions and notation of the previous section. 
For the sake of notational simplicity, in the following we shall assume that 
$\E[f(X_1)]=0$. Also, we shall assume that 
$X_1\sim N_d(0, I_d)$, implying that
$X^{(1)}_k, \ldots, X^{(d)}_k$ are independent $N(0,1)$ random variables for all $k\in \Z$.
Note that this 
last assumption is not restrictive:
indeed, by reduction
of variables and at the cost of possibly decreasing the value of $d$, we may always assume that $X_k\sim N_{d}(0, \Sigma )$ for an {\sl invertible} matrix $\Sigma \in \R^{d \times d}$,
and by a linear transformation we can further restrict ourselves to the case $X_k\sim N_{d}(0, I_{d})$. Observe that, if $f(x)$ has Hermite rank $q$ with respect to $X_1$, then $f(Lx)$ has the same Hermite rank $q$ with respect to $L^{-1}X_1$, for every invertible matrix $L$ -- see also \cite[p. 2249]{AR}.

Now, let $\Lambda$ denote the set of all vectors
$\alpha=(\alpha_1,  \ldots, \alpha_d)$ with $\alpha_i\in \N \cup \{0\}$. For any multi-index
$\alpha \in \Lambda $, we introduce
the notation $|\alpha | = \sum_{i= 1}^d \alpha _i$ and $\alpha ! = \prod_{i= 1}^d \alpha _i!$. When $\E[f^2(X_1)]<\infty $,
the function $f$ possesses the unique {\sl Hermite expansion}
\begin{equation} \label{hermite}
f(x)= \sum_{\alpha \in \Lambda} a_{\alpha} \prod_{i=1}^d H_{\alpha_i} (x_i), \qquad
a_{\alpha} = (\alpha !)^{-1} \E\Big[f(X_1) \prod_{i=1}^d H_{\alpha_i} (X_1^{(i)})\Big],
\end{equation}
where $(H_j)_{j\geq 0}$ is the sequence of Hermite polynomials, recursively defined as: $H_0 =1 $, and
\[
H_j = \delta  H_{j-1}, \quad j\geq 1,
\]
where $\delta f(x) = xf(x) - f'(x)$ (for instance: $H_1(x) = x$, $H_2(x) = x^2 -1$, and so on).

The following well-known statement provides a further characterization of Hermite ranks.

\begin{prop}\label{D : HermiteRank} Let the notation and assumptions of this section prevail 
(in particular, 
$\E[f(X_1)]=0$ and 
$X_1 \sim N_d(0, I_d)$).
Then, the function $f$
has Hermite rank $q\geq 1$ if and only if $a_{\alpha}=0$ for all
$\alpha \in \Lambda $ with $|\alpha |<q$ and $a_{\alpha}\not =0$ for some $\alpha \in \Lambda $ with $|\alpha |=q$.
In particular, if $f$ has Hermite rank $q$, then its Hermite expansion has the form
\begin{equation} \label{fm}
f(x)= \sum_{m=q}^\infty f_m(x), \qquad f_m(x)= \sum_{\alpha \in \Lambda: |\alpha|=m} a_{\alpha} \prod_{i=1}^d H_{\alpha_i} (x_i).
\end{equation}
\end{prop}

\noindent {\bf Remark on notation.} Fix a function $f$ such that $\E[f(X_1^2)]<\infty$ and $f$ has Hermite rank $q\geq1$. Our main results are expressed in terms of the following collection of coefficients (\ref{thetaJ})--(\ref{gamma}):

\begin{eqnarray}
\label{thetaJ}\theta(j)\!\!&=&\!\! \max_{1\leq i,l\leq d} |r^{(i,l)}(j)| \\[1.5 ex]
\label{K} K\!\!&=&\!\!\inf_{k\in \N} \{ \theta(j) \leq d^{-1}, \forall |j|\geq k\} \quad (\mbox{with } \inf\emptyset=\infty),\\[1.5 ex]
\label{theta}\theta\!\!&=&\!\!  \sum_{j\in \Z} \theta(j)^q, \\[1.5 ex]
\label{sigmam}\sigma_m^2\!\! &=&\!\!  \E[f^2_m(X_1)] +2\sum_{k=1}^\infty \E[f_m(X_1)f_m(X_{1+k})] \quad (\mbox{for }m\geq q),  \\[1.5 ex]
\label{gamma}  \gamma_{n,m,e}\!\! &=&\!\! \sqrt{2\theta n^{-1}
 \sum_{|j|\leq n} \theta(j)^e \sum_{|j|\leq n} \theta(j)^{m-e}} \,\,\, (\mbox{for }m\geq q \mbox{ and } 1\leq e\leq m-1).
\end{eqnarray}
The coefficients $\theta(j)$, $K$, $\theta$, $\sigma^2_m$ and $\gamma_{n,m,e}$ will be also combined into the following expressions (\ref{A1})--(\ref{A5}):
\begin{eqnarray}
\label{A1}A_{1,n}\!\!&=&\!\! \frac {\E[f^2(X_1)] }{2} \left[ \frac{2K^2}{n}+
d^q\left(\sum_{|j|\leq n} \theta(j)^q \frac{|j|}{n} +\sum_{|j|>n}\theta(j)^q\right)
\right]\\[1.5 ex]
\label{A2} A_{2,N}\!\!&=& \!\! 2(2K+d^q \theta)
\sqrt{\E[f^2(X_1)]\sum_{m=N+1}^\infty \E[f^2_m(X_1)]} \\[1.5 ex]
\label{A3} A_{3,n,N}\!\!&=&\!\! \frac{\E[f^2(X_1)]}{2}\sum_{m=q}^N \left(\frac{d^m }{m m!} \sum_{l=1}^{m-1}ll!\binom{m}{l}^2\sqrt{(2m-2l)!}
\,\gamma_{n,m,l}\right) \\[1.5 ex]
\label{A4} A_{4,n,N}\!\!&=&\!\! \frac{\E[f^2(X_1)](2K+d^q \theta)^{1/2}}{2}\times \\
&&\quad\quad\quad\quad\quad\quad\times\sum_{q\leq p<s\leq N} d^{s/2}\sqrt{\frac{p!}{s!}} \frac{p+s}{p}\binom{s-1}{p-1}\sqrt{(s-p)!}\notag
\gamma_{n,s,s-p}^{1/2} \\[1.5 ex]
 A_{5,n,N}\!\!&=&\!\!\frac{\E[f^2(X_1)]}{2\sqrt{2}} \sum_{q\leq p<s\leq N} (p+s)
\sum_{l=1}^{p-1} (l-1)!\binom{p-1}{l-1}\binom{s-1}{l-1}\sqrt{(p+s-2l)!}\times \nonumber \\[1.5 ex]
&\times & \left(\frac{d^s}{s!}
\gamma_{n,s,s-l}+\frac{d^p}{p!} \gamma_{n,p,p-l}\right). \label{A5}
\end{eqnarray}

Note that the coefficients $K$, $\theta$ and $\sigma_m$ can in general be infinite, and also that, if $\E[f(X_1)^2]<\infty$, if $f$ has Hermite rank $q$, and if (\ref{zoe}) is in order, then
\begin{equation}\label{tMOE}
\sigma^2 = \sum_{m=q}^\infty \sigma^2_m<\infty,
\end{equation}
where $\sigma^2$ is defined in Theorem \ref{th1}.

The next statement, which is the main result of the paper, asserts that the quantities defined above can be used to write explicit bounds of the type (\ref{TEST1}).

\begin{theo}[Quantitative Breuer-Major Theorem] \label{mainth}
Let the notation and assumptions of this section prevail 
(in particular, $\E[f(X_1)]=0$ and $X_1 \sim N_d(0, I_d)$), and
assume that the conditions of Theorem \ref{th1} are satisfied. 
Then, the coefficients appearing in formulae (\ref{thetaJ})--(\ref{A5}) are all finite. 
Moreover, the following three bounds are in order.\\
\\
{\bf (1)} For any function
$h\in C^2(\R)$ (that is, $h$ is twice continuously differentiable) with bounded second derivative, and for every $n>K$,
\begin{equation} \label{mainestimate}
|\E[h(S_n) ] -\E[h(S)]|\leq \|h''\|_\infty\Big(A_{1,n}+\inf_{N\geq q}\{A_{2,N}+A_{3,n,N}+A_{4,n,N}+A_{5,n,N}\}\Big).
\end{equation}
{\bf (2)} For any Lipschitz function $h$, and for every $n>K$,
\begin{eqnarray} \label{mainestimate2}
&&|\E[h(S_n) ] -\E[ h(S)]|\leq \frac{\|h'\|_\infty}{2} \times\\
&& \times \left\{ A_{2,n}\left[
\frac{1}{\sigma }+\frac{1}{\sqrt{\left( 2K+d^{q}\theta \right) \mathbb{E}%
\left[ f^{2}\left( X_{1}\right) \right] }}\right] +4\inf_{N\geq q}\frac{%
A_{1,n}+A_{3,n,N}+A_{4,n,N}+A_{5,n,N}}{\sqrt{\sum_{m=q}^{N}\sigma _{m}^{2}}}%
\right\}.\notag
\end{eqnarray}
{\bf (3)} For any $z\in\R$, and for every $n>K$,
\begin{eqnarray} \label{mainestimate3}
&&|\mathbb{P}(S_n\leq z) ] -\mathbb{P}(S\leq z)|\leq \\
&&\frac{\sqrt{2}}{\sigma}\sqrt{A_{2,n}\left[
\frac{1}{\sigma }+\frac{1}{\sqrt{\left( 2K+d^{q}\theta \right) \mathbb{E}%
\left[ f^{2}\left( X_{1}\right) \right] }}\right] +4\inf_{N\geq q}\frac{%
A_{1,n}+A_{3,n,N}+A_{4,n,N}+A_{5,n,N}}{\sqrt{\sum_{m=q}^{N}\sigma _{m}^{2}}}}.\notag
\end{eqnarray}
\end{theo}

We will now demonstrate that Theorem \ref{mainth} implies a stronger version of Theorem \ref{th1}, namely that the convergence (\ref{clt}) takes place with respect to topologies that are stronger than the one of convergence in distribution. To prove this claim,
we need to show in particular that, under the assumptions of Theorem \ref{th1}, $\gamma_{n,m,e} \rightarrow 0$ as $n\rightarrow \infty$ for any choice of $m\geq q$ and $1\leq e\leq m-1$.
This is a consequence of the next Lemma \ref{lemconv}. In what follows, given positive sequences $b_n,c_n$, $n\geq 1$, we shall write $b_n \lesssim c_n$ whenever $b_n/c_n$ is bounded, and
$b_n \sim  c_n$ if $b_n \lesssim c_n$ and $c_n \lesssim b_n$.

\begin{lem} \label{lemconv}
Let $(a_k)_{k \in \Z}$ be a sequence of positive real numbers such that $\sum_{k\in \Z} a_k^m<\infty$ for some
$m\in \N$. If $1\leq e\leq m-1$, then
\begin{equation*}
n^{-1+\frac em} \sum_{|k|\leq n} a_k^e \rightarrow 0.
\end{equation*}
\end{lem}
\textit{Proof.}
Fix $\delta \in (0,1)$, and decompose the sum as $\sum_{k=1}^n = \sum_{k=1}^{[n\delta]} + \sum_{k=[n\delta]+1}^n$.
By the H\"{o}lder inequality we obtain (recall that $\sum_{k=1}^\infty a_k^m$ is finite)
\[
n^{-1+e/m} \sum_{k=1}^{[n\delta]} a_k^e \leq  n^{-1+e/m} \left(n\delta\right)^{1-e/m} \left(\sum_{k=1}^\infty a_k^m\right)^{e/m} \leq c\delta^{1-e/m},
\]
where $c$ is some constant, as well as
\[
n^{-1+e/m} \sum_{k=[n\delta]+1}^n a_k^e \leq \left(\sum_{k=[n\delta]+1}^n a_k^m\right)^{e/m}.
\]
The first term converges to 0 as $\delta$ goes to zero  (because $1\leq e\leq m-1$), and the second also converges to 0 for fixed $\delta$ and $n\rightarrow\infty$. This proves the claim.
\qed

Now recall that, if $X,Y$ are two real-valued random variables, then the {\sl Kolmogorov distance} between the law of $X$ and the law of $Y$ is given by
\begin{eqnarray}
\label{KOL} d_{Kol}(X,Y) = \sup_{z\in\R}|\PP(X\leq z) - \PP(Y\leq z)|.
\end{eqnarray}
If $\E|X|,\E|Y|<\infty$, one can also meaningfully define the {\sl Wasserstein distance}
\begin{eqnarray}
\label{WASS} d_{W}(X,Y) = \sup_{f\in {\rm Lip}(1)}|\E[f(X)] - \E[f(Y)]|,
\end{eqnarray}
where ${\rm Lip}(1)$ indicates the collection of all Lipschitz functions with Lipschitz constant $ \leq 1$. Finally, if $X,Y$ have finite second moments, for every constant $C>0$ one can define the distance
\begin{eqnarray}
\label{WASS2} d_{C}(X,Y) = \sup_{f\in \mathscr{D}^2_C}|\E[f(X)] - \E[f(Y)]|,
\end{eqnarray}
where $\mathscr{D}^2_C$ stands for the class of all twice continuously differentiable functions having a second derivative bounded by $C$. Note that the topologies induced by $d_{Kol}$, $d_{W}$ and $d_{C}$, on the probability measures on $\R$, are strictly stronger than the topology of convergence in distribution (see e.g. \cite[Ch. 11]{Dudley book}).

The next consequence of Theorem \ref{mainth} provides the announced refinement of Theorem \ref{th1}.

\begin{cor}[Breuer-Major, strong version]\label{C : topologies} 
Let the notation and assumptions of this section prevail 
(in particular, $\E[f(X_1)]=0$ and $X_1 \sim N_d(0, I_d)$), and
assume that the conditions of Theorem \ref{th1} are satisfied. 
Then, the convergence (\ref{clt}) takes place with respect to the three distances $d_{Kol}$, $d_{W}$ and $d_{C}$ (for all $C>0$), namely
\[
\lim_{n\rightarrow\infty}d_{Kol}(S_n,S)=\lim_{n\rightarrow\infty} d_{W}(S_n,S)=\lim_{n\rightarrow\infty}d_{C}(S_n,S) = 0.
\]
\end{cor}
{\it Proof.}
Under the assumptions of Theorem \ref{th1}, one has that $A_{1,n}\to 0$ as $n\to\infty$ (because $\theta<\infty$, $\sum_{|j|\leq n}\theta(j)^q\frac{|j|}{n}\to 0$
as $n\to\infty$ by bounded convergence). On the other hand, because of (\ref{tMOE}) and since $\E[f^2(X_1)]<\infty $, one has that
$A_{2,N}\rightarrow 0$ as $N\rightarrow \infty $. Moreover, since
$\gamma_{n,m,e} \rightarrow 0$ for any $m\geq q$ and $1\leq e\leq m-1$  (due to Lemma \ref{lemconv}), one has
that $A_{j,n,N}\rightarrow 0$, $j=3,4,5$,
for any fixed $N$ as $n\rightarrow \infty$. We deduce that $\inf_{N\geq q}\{A_{2,N}+A_{3,n,N}+A_{4,n,N}+A_{5,n,N}\}\to 0$ as $n\to\infty$.
To conclude the proof, it remains to apply (\ref{mainestimate})--(\ref{mainestimate3}).
\qed

\medskip

Next, we present a simplified version of Theorem \ref{mainth} for  $d=1$ and $f=H_q$, where $H_q$
is the $q$th Hermite polynomial. Notice that in this case $K=0$.

\begin{cor} [Hermite subordination] \label{cor1}
Assume that $d=1$, $f=H_q$ and  $\sum_{j \in \Z} |r(j)|^q<\infty$.
\begin{description}
\item[\bf (1)] For any function
$h\in C^2(\R)$ with bounded second derivative it holds that
\begin{equation} \label{mainestimatecor}
|\E[h(S_n) ] -\E[ h(S)]|\leq \|h''\|_\infty\Big(A_{1,n}+A_{3,n}\Big).
\end{equation}
\item[\bf (2)] For any Lipschitz function $h$ it holds that
\begin{eqnarray} \label{mainestimate2cor}
|\E[h(S_n) ] -\E[ h(S)]|\leq \frac{2\|h'\|_\infty}{\sigma} \Big(A_{1,n}+A_{3,n}\Big).
\end{eqnarray}
\item[\bf (3)] For any $z\in \R$ it holds that
\begin{eqnarray} \label{mainestimate3cor}
|\mathbb{P}(S_n\leq z) ] -\mathbb{P}(S\leq z)|\leq \frac{2\|s_z'\|_\infty}{\sigma} \Big(A_{1,n}+A_{3,n}\Big),
\end{eqnarray}
where $s_z$ is the solution of the Stein's equation associated with the function $h(x)=1_{(-\infty, z]}(x)$,
i.e. $s_z$ solves the differential equation
$$ 1_{(-\infty, z]}(x) - \Phi (z)=s_z'(x) - xs_z(x), \qquad x\in \R,$$
with $\Phi$ being the distribution function of $N(0,1)$. Furthermore, we have that $\|s_z'\|_\infty\leq 1$ for all $z\in \R$.
\end{description}
In this context, the constants $A_{1,n}$ and $A_{3,n}$ are given by
\begin{eqnarray*}
A_{1,n}&=& \frac {q!}{2} \theta \left(
\sum_{|j|\leq n} |r(j)|^q \frac{|j|}{n} +\sum_{|j|>n}|r(j)|^q
\right),\nonumber \\[1.5 ex]
A_{3,n}&=& \frac{1}{2q} \sum_{l=1}^{q-1}ll!\binom{q}{l}^2\sqrt{(2q-2l)!}
\,\gamma_{n,q,l},
\end{eqnarray*}
with $\gamma_{n,q,l}$ defined by (\ref{gamma}).
\end{cor}
\textit{Proof.} From Theorem 3.1 in \cite{NP1} and Theorem \ref{ProbApp} in Section 3.3  we obtain the estimate
\begin{equation*}
|\E[h(S_n) ] -\E[ h(S)]| \leq c_h\, \E|\sigma^2- \langle DF,-DL^{-1}F \rangle_\HH |,
\end{equation*}
where $c_h=\frac{\| h''\|_\infty}{2}$ in (1), $c_h=\frac{\| h'\|_\infty}{\sigma}$ in (2) and $c_h=\frac{\| s_z'\|_\infty}{\sigma}$
in (3). We readily deduce the assertion since $\E|\sigma^2- \langle DF,-DL^{-1}F \rangle_\HH |\leq 2(A_{1,n}+A_{3,n})$,
which follows from the proof of Theorem \ref{mainth}. \qed \\ \\
We remark that the upper bound in (\ref{mainestimate3cor}) of Corollary \ref{cor1} is stronger than the general upper bound obtained
in  (\ref{mainestimate3}) of Theorem \ref{mainth}.

Next, we apply Corollary \ref{cor1} to some particular classes of covariance functions $r$.

\begin{ex}  \label{fBM} \rm \textbf{(Covariance functions with polynomial decay)}
Assume that $d=1$ and $f=H_q$ with $q\geq 2$, and consider a covariance function $r$ which is regular varying with parameter $a<0$.
That is, for all $|k|\geq 1$, $|r(k)| = |k|^a l(|k|)$, where $l$ is a slowly varying function. Recall that for any regular varying function
$r$ with parameter $\alpha<0$, we have the following discrete version of \textit{Karamata's theorem} (see e.g. \cite{BGT}):
\begin{eqnarray*}
\alpha>-1&:& \qquad \frac{\sum_{k=1}^n |r(k)|}{n^{a+1} l(n)} \rightarrow 1/(\alpha+1) , \\[1.5 ex]
\alpha<-1&:& \qquad  \frac{\sum_{k=n}^\infty  |r(k)|}{n^{a+1} l(n)} \rightarrow -1/(\alpha+1) ,
\end{eqnarray*}
as $n\rightarrow \infty $. Assume now that $a<-\frac{1}{q}$, which implies that the conditions of Theorem \ref{mainth}
are satisfied, and $ae \not=-1$ for any $e=1,\ldots,q-1$. By the afore-mentioned convergence results we immediately deduce the following
estimates ($1\leq e\leq q-1$)
\begin{eqnarray*}
\sum_{|j|\leq n} |r(j)|^q \frac{|j|}{n} &\lesssim& n^{-1}+ n^{aq+1}l(n), \\
\sum_{|j|>n}|r(j)|^q &\lesssim& n^{aq+1} l(n), \\
\gamma_{n,q,e} &\lesssim & n^{-1/2}+ n^{a/2}l(n).
\end{eqnarray*}
Thus, for all three cases of Corollary \ref{cor1} we conclude that
\begin{equation*}
|\E[h(S_n) ] -\E[ h(S)]| \lesssim \left \{ \begin{array} {cc}
n^{-1/2}: & a< -1 \\[1.5 ex]
n^{a/2}l(n): & a \in (-1, -\frac{1}{q-1}) \\[1.5 ex]
n^{\frac{aq+1}{2}}l(n): & a \in ( -\frac{1}{q-1}, -\frac{1}{q})
\end{array} \right.
\end{equation*}
Clearly, the same estimates hold for $d(S_n, S)$, where $d=d_{Kol}$, $d=d_W$ or $d=d_C$.
\end{ex}

\begin{ex} \rm \textbf{(The fractional Brownian motion case)}
Let $d=1$, $f=H_q$ with $q\geq 2$ and consider the fractional Gaussian noise $X_i=B_i^H-B_{i-1}^H$, where $B^H$ is a fractional Brownian
motion with parameter $H\in (0,1)$. Recall that $B^H=(B_t^H)_{t\geq 0}$ is a centered Gaussian process (with stationary increments) with covariance structure given by
$$\mathbb E[B_t B_s]= \frac 12 (|t|^{2H} + |s|^{2H} - |t+s|^{2H})$$
It is well-known that the correlation function $r$ of the fractional Brownian noise has the following
form:
\begin{equation*}
|r(k)| = |k|^{2H-2} l(|k|), \qquad k\geq 1,
\end{equation*}
with $l(|k|) \rightarrow 2H |2H-1|$ as $|k|\rightarrow \infty$ when $H \not= \frac 12$, and $l(|k|) =0$ for $|k|\geq 1$
when $H = \frac 12$.
As in the previous example we immediately deduce that
\begin{equation*}
|\E[h(S_n) ] -\E[ h(S)]| \lesssim \left \{ \begin{array} {cc}
n^{-1/2}: & a\in (-2, -1] \\[1.5 ex]
n^{a/2}: & a \in [-1, -\frac{1}{q-1}] \\[1.5 ex]
n^{\frac{aq+1}{2}}: & a \in [ -\frac{1}{q-1}, -\frac{1}{q})
\end{array} \right.
\end{equation*}
with $a=2H-2$ and the same estimates hold for $d(S_n, S)$ with $d=d_{Kol}$, $d=d_W$ or $d=d_C$.
Let us remark that these upper bounds coincide with those derived in Theorem 4.1 in \cite{NP1}.

We finally remark that the rate $n^{-1/2}$ for $a=2H-2\in (-2, -1]$ has been proved to be optimal in \cite{NP2}.
For the other two cases the optimality question is still an open problem.
\end{ex}

\section{Toolbox}\label{S : toolbox}
\subsection{Malliavin calculus on a Gaussian space}
We shall now provide a short introduction to the tools of Malliavin
calculus that are needed in the proof of our main Theorem \ref{mainth}. The
reader is referred to \cite{N} for any unexplained definition or result. Let $\HH$ be a real separable Hilbert space. We denote by $W =
\{W(h):h \in \HH\}$ an \textsl{isonormal Gaussian process} over
$\HH$, that is, $W$ is a centered Gaussian family indexed by
the elements of $\HH$ and such that, for every $g_1,g_2\in\HH$,
\begin{equation}\label{isometry}
\E\big[W(g_1)W(g_2)\big]=\langle g_1,g_2\rangle_\HH.
\end{equation}
In what follows, we shall use the notation ${L}^2(W)$ $=$
$L^2(\Omega,\sigma(W),\PP)$. For every $q\geq 1$, we write
$\HH^{\otimes q}$ to indicate the $q$th tensor power of $\HH$; the
symbol $\HH^{\odot q}$ indicates the $q$th \textsl{symmetric}
tensor power of $\HH$, equipped with the norm
$\sqrt{q!}\|\cdot\|_{\HH^{\otimes q}}$. We denote by $I_q$ the
isometry between $\HH^{\odot q}$ and the $q$th Wiener chaos of
$X$. It is well-known (see again \cite[Ch. 1]{N}) that any random variable $F$ belonging to $L^2(W)$ admits the \textsl{chaotic expansion}:
\begin{equation}\label{ChaosExpansion}
F=\sum_{q=0}^\infty I_q(f_q),
\end{equation}
where $I_0(f_0):=E[F]$, the series converges in $L^2$ and
the kernels $f_q\in\HH^{\odot q}$, $q\geq1$, are uniquely
determined by $F$. In the particular case where $\HH=L^2(A,\mathscr{A},\mu)$, with $(A,\mathscr{A})$ a measurable
space and $\mu$ a $\sigma$-finite and non-atomic measure, one
has that $\HH^{\odot q}= L^2_s(A^q,\mathscr{A}^{\otimes q},
\mu^{\otimes q})$ is the space of symmetric and square integrable
functions on $A^q$. Moreover, for every $f\in\HH^{\odot q}$,
$I_q(f)$ coincides with the multiple Wiener-It\^o integral (of
order $q$) of $f$ with respect to $W$ (see \cite[Ch. 1]{N}).
It is well-known that a random variable of the type $I_q(f)$,
$f\in\HH^{\odot q}$, has finite moments of all orders (see
\cite[Ch. VI]{Janson}). For every
$q\geq 0$, we write $J_q$ to indicate the orthogonal projection
operator on the $q$th Wiener chaos associated with $W$, so that,
if $F \in L^2(W)$ is as in (\ref{ChaosExpansion}),
then $J_q F=I_q(f_q)$ for every $q\geq 0$.

Let $\{e_k,\,k\ge 1\}$ be a complete orthonormal system in $\HH$.
Given $f\in \HH^{\odot p}$ and $g\in\HH^{\odot q}$, for every
$r=0,\ldots,p\wedge q$, the $r$th contraction of $f$ and $g$ is
the element of $\HH^{\otimes (p+q-2r)}$ defined as
\begin{equation}\label{defContr}
f\otimes_r g = \sum_{i_1,\ldots,i_r=1}^\infty \langle
f,e_{i_1}\otimes \ldots\otimes e_{i_r}\rangle_{\HH^{\otimes r}}
\otimes \langle g,e_{i_1}\otimes \ldots\otimes
e_{i_r}\rangle_{\HH^{\otimes r}}.
\end{equation}
In the particular case where $\HH=L^2(A,\mathscr{A},\mu)$
(with $\mu$ non-atomic), one has that
$$
f\otimes_r g = \int_{A^r}
f(t_1,\ldots,t_{p-r},s_1,\ldots,s_r)\,g(t_{p-r+1},\ldots,
t_{p+q-2r},s_1,\ldots,s_r)d\mu(s_1)\ldots d\mu(s_r).
$$
Moreover, $f\otimes_0 g=f\otimes g$ equals the tensor product of
$f$ and $g$ while, for $p=q$, $f\otimes_p g=\langle
f,g\rangle_{\HH^{\otimes p}}$. Note that, in general, the contraction $f\otimes_r g$ is \textsl{not}
a symmetric element of $\HH^{\otimes (p+q-2r)}$. The canonical
symmetrization of $f\otimes_r g$ is written
$f\widetilde{\otimes}_r g$. The following
multiplication formula is also very useful: if $f\in \HH^{\odot p}$ and
$g\in\HH^{\odot q}$, then
\begin{equation}\label{multiplication}
I_p(f)I_q(g)=\sum_{r=0}^{p\wedge q} r!\binom{p}{r}\binom{q}r
I_{p+q-2r}(f\widetilde{\otimes}_r g).
\end{equation}
Let $\mathscr{S}$ be the set of all smooth cylindrical random
variables of the form $$F = g\big(W(\phi_1), \ldots,
W(\phi_n)\big),$$ where $n\geq 1$, $g : \R^n \rightarrow \R$ is a
smooth function with compact support and $\phi_i\in\EuFrak H$. The
Malliavin derivative of $F$ with respect to $W$ is the element of
$L^2(\Omega, \EuFrak H)$ defined as
$$DF\; =\; \sum_{i =1}^n \frac{\partial g}{\partial x_i}\big(W(\phi_1), \ldots, W(\phi_n)\big)
\phi_i.$$ Also, $DW(\phi) = \phi$ for every $\phi\in \HH$.
As usual, $\sk^{1,2}$ denotes the
closure of $\mathscr{S}$ with respect to the norm $\| \cdot
\|_{1,2}$, defined by the relation
$$\| F\|_{1,2}^2 \; = \; \E\lcr F^2\rcr +
\E\big[ \| D F\|_{\HH}^2\big].$$
Note that, if $F$ is equal to a finite sum of multiple Wiener-Itô
integrals, then $F\in\sk^{1,2}$. The Malliavin
derivative $D$ verifies the following \textsl{chain rule}: if
$\varphi:\R^n\rightarrow\R$ is in $\mathscr{C}^1_b$ (that is, the
collection of continuously differentiable functions with
bounded partial derivatives) and if $\{F_i\}_{i=1,\ldots,n}$ is a vector of
elements of $\sk^{1,2}$, then
$\varphi(F_1,\ldots,F_n)\in\sk^{1,2}$ and
$$
D\varphi(F_1,\ldots,F_n)=\sum_{i=1}^n
\frac{\partial\varphi}{\partial x_i} (F_1,\ldots, F_n)DF_i.
$$
We denote by $\delta$ the
adjoint of the operator $D$, also called the \textsl{divergence
operator}. A random element $u \in L^{2}(\Omega, \HH)$
belongs to the domain of $\delta$, noted ${\rm Dom}\delta$, if
and only if it verifies
$$
| \E \langle D F,u\rangle_{\EuFrak H}|\leq c_u\,\|F\|_{L^2}\quad\mbox{for any }F\in{\mathscr{S}},
$$
where $c_u$ is a constant depending only on $u$. If $u \in
{\rm Dom}\delta $, then the random variable $ \delta(u)$ is
defined by the duality relationship (sometimes called `integration by
parts formula'):
\begin{equation}\label{ipp}
\E [F \delta(u)]=  \E \langle D F, u \rangle_{\HH},
\end{equation}
which holds for every $F \in \sk^{1,2}$.

The operator $L$, acting on square integrable random variables of
the type (\ref{ChaosExpansion}), is defined through the projection
operators $\{J_q\}_{q \geq 0}$ as $L=\sum_{q=0}^\infty -qJ_q,$ and
is called the \textsl{infinitesimal generator of the
Ornstein-Uhlenbeck semigroup}. It verifies the following crucial
property: a random variable $F$ is an element of $\rm{Dom}$$ L
\,\,(=\sk^{2,2})$ if, and only if, $F\in{\rm Dom}\delta D$ (i.e.
$F\in\sk^{1,2}$ and $DF\in{\rm Dom}\delta$), and in this case:
$\delta DF=-LF.$  Note that a random variable $F$ as in
(\ref{ChaosExpansion}) is in $\sk^{1,2}$  if
and only if
$$
\sum_{q=1}^\infty (q+1)!\|f_q\|^2_{\HH^{\otimes q}}<\infty,
$$
and also $\E\big[\|DF\|^2_\HH\big]=\sum_{q\geq 1}
qq!\|f_q\|^2_{\HH^{\otimes q}}$. If $\HH=L^2(A,\mathscr{A},\mu)$
(with $\mu$ non-atomic), then the derivative of a random variable
$F$ as in (\ref{ChaosExpansion}) can be identified with the
element of $L^2(A\times \Omega)$ given by
\begin{equation}\label{dtf}
D_a F=\sum_{q=1}^\infty qI_{q-1}\big(f_q(\cdot,a)\big), \quad a
\in A .
\end{equation}
We also define the operator $L^{-1}$, which is the
\textsl{pseudo-inverse} of $L$, as follows: for every $F\in{ L^2}(W)$,
we set $L^{-1}F$ $=$ $\sum_{q\geq 1}\frac{1}{q}
J_q(F)$. Note that $L^{-1}$ is an operator with values in
$\sk^{2,2}$ and that $LL^{-1}F=F-\E[F]$ for all $F\in L^2(W)$.

\subsection{Assessing norms and scalar products}

The following statement plays a crucial role in the proof of Theorem \ref{mainth}.

\begin{lem} \label{lem1}
Let $F=I_p(h)$ and $G=I_s(g)$ with $h\in \HH^{\odot p}$, $g\in \HH^{\odot s}$ and $p<s$ ($p,s\geq 1$).
Then
\begin{equation}\label{varident1}
{\rm Var}\left[\frac1s\|DG\|^2_{\HH}\right]=
\frac{1}{s^2}\sum_{l=1}^{s-1}l^2l!^2\binom{s}{l}^4 (2s-2l)!\|g\widetilde{\otimes}_{l} g\|^2_{\HH^{\otimes 2s-2l}},
\end{equation}
and
\begin{eqnarray}
&&\E\left[\left(\frac1s\left\langle DF,DG\right\rangle_{\HH}
\right)^2\right]
\leq p!\binom{s-1}{p-1}^2(s-p)!\E[F^2]\|g\otimes_{s-p}g\|_{\HH^{\otimes 2p}}\label{varident2}\\[1.5 ex]
&&+
\frac{p^2}{2}\sum_{l=1}^{p-1}(l-1)!^2\binom{p-1}{l-1}^2\binom{s-1}{l-1}^2(p+s-2l)!\Big(
\|h\otimes_{p-l}h\|^2_{\HH^{\otimes
2l}}+\|g\otimes_{s-l}g\|^2_{\HH^{\otimes 2l}}\Big).\nonumber
\end{eqnarray}
\end{lem}
\textit{Proof.} [Proof of (\ref{varident1})] We have
$DG=sI_{s-1}(g)$ so that, by using (\ref{multiplication})
\begin{eqnarray*}
\frac1s\|DG\|^2_{\HH}&=&s\|I_{s-1}(g)\|^2_{\HH}
=s\sum_{l=0}^{s-1}l!\binom{s-1}{l}^2 I_{2s-2-2l}(g\widetilde{\otimes}_{l+1} g)\\[1.5 ex]
&=&s\sum_{l=1}^{s}(l-1)!\binom{s-1}{l-1}^2 I_{2s-2l}(g\widetilde{\otimes}_{l} g)\\[1.5 ex]
&=&s!\|g\|^2_{\HH^{\otimes s}}+s\sum_{l=1}^{s-1}(l-1)!\binom{s-1}{l-1}^2 I_{2s-2l}(g\widetilde{\otimes}_{l} g)\\[1.5 ex]
&=&\E[G^2]+s\sum_{l=1}^{s-1}(l-1)!\binom{s-1}{l-1}^2 I_{2s-2l}(g\widetilde{\otimes}_{l} g).
\end{eqnarray*}
The orthogonality property  of multiple integrals leads to (\ref{varident1}).

\noindent[Proof of (\ref{varident2})] Thanks once again to (\ref{multiplication}), we can write
\begin{eqnarray*}
\langle DF,DG\rangle_{\HH} &=&p\,s\left\langle I_{p-1}(h),I_{s-1}(g)\right\rangle_{\HH}\\
&=&p\,s\sum_{l=0}^{p\wedge s-1} l!\binom{p-1}{l}\binom{s-1}{l} I_{p+s-2-2l}(h\widetilde{\otimes}_{l+1}g)\\
&=&p\,s \sum_{l=1}^{p\wedge s} (l-1)!\binom{p-1}{l-1}\binom{s-1}{l-1} I_{p+s-2l}(h\widetilde{\otimes}_l g).
\end{eqnarray*}
It follows that
\begin{eqnarray}
\E\left[\left(\frac1s\left\langle DF,DG\right\rangle_{\HH}\right)^2\right] \notag
=p^2\sum_{l=1}^{p}(l-1)!^2
\binom{p-1}{l-1}^2\binom{s-1}{l-1}^2 (p+s-2l)!
\|h\widetilde{\otimes}_l g\|^2_{\HH^{\otimes (p+s-2l)}} .\\
\label{Murray}
\end{eqnarray}
If $l<p$, then
\begin{eqnarray*}
\|h\widetilde{\otimes}_l g\|^2_{\HH^{\otimes (p+s-2l)}}
&\leq& \|h\otimes_l g\|^2_{\HH^{\otimes (p+s-2l)}}
=\langle h\otimes_{p-l} h, g\otimes_{s-l}g\rangle_{\HH^{\otimes 2l}}\\
&\leq&
\|h\otimes_{p-l}h\|_{\HH^{\otimes 2l}}\|g\otimes_{s-l}g\|_{\HH^{\otimes 2l}}\\
&\leq&\frac12\left(
\|h\otimes_{p-l}h\|_{\HH^{\otimes 2l}}^2+\|g\otimes_{s-l}g\|_{\HH^{\otimes 2l}}^2
\right).
\end{eqnarray*}
If $l=p$, then
$$
\|h\widetilde{\otimes}_p\, g\|^2_{\HH^{\otimes (s-p)}} \leq
\|h\otimes_p \,g\|^2_{\HH^{\otimes (s-p)}} \leq
\|h\|^2_{\HH^{\otimes p}}\|g\otimes_{s-p}g\|_{\HH^{\otimes 2p}}.
$$
By plugging these last expressions into (\ref{Murray}), we deduce
immediately (\ref{varident2}).
\qed
\subsection{Estimates via interpolations and Stein's method}

The forthcoming Theorem \ref{ProbApp} contains two bounds on normal approximations, that are expressed in terms of Malliavin operators. As anticipated, the proof of Point (1) uses an interpolation technique already applied in \cite{NouPeReinertAOP, PecZheng}, which is close to the `smart path method' of Spin Glasses \cite{talag}. Point (2) uses estimates from \cite{NP1}.

\begin{theo}\label{ProbApp}
Let $F$ be a centered element of $\sk^{1,2}$ and let $Z\sim N(0,\sigma^2)$, $\sigma>0$.
\begin{description}
\item[\bf (1)] Suppose that $h : \R\rightarrow\R$ is twice continuously differentiable and has a bounded second derivative. Then,
\begin{equation}\label{E : est1}
\big| \E[h(F)]- \E[h(Z)] \big| \leq \frac{\| h''\|_\infty}{2}\, \E|\sigma^2- \langle DF,-DL^{-1}F \rangle_\HH |.
\end{equation}
\item[\bf (2)] Suppose that $h : \R\rightarrow\R$ is Lipschitz. Then,
\begin{equation}\label{E : est2}
\big| \E[h(F)]- \E[h(Z)] \big| \leq \frac{\| h'\|_\infty}{\sigma}\, \E|\sigma^2- \langle DF,-DL^{-1}F \rangle_\HH |.
\end{equation}
\end{description}
\end{theo}
{\it Proof.} {\bf (1)} Without loss of generality, we may assume that $F$ and $Z$ are independent and defined on the same probability space. Fix $h$ as in the statement, and define the function $\Psi(t) = \E[h(\sqrt{1-t}F+\sqrt{t}Z)]$, $t\in [0,1]$. Standard results imply that $\Psi$ is differentiable for every $t\in (0,1)$, and that
\[
\Psi'(t) = \frac{1}{2\sqrt{t}}\, \E[h'(\sqrt{1-t}F+\sqrt{t}Z)Z] - \frac{1}{2\sqrt{1-t}} \, \E[h'(\sqrt{1-t}F+\sqrt{t}Z)F].
\]
By using independence and integration by parts, we obtain immediately that
\[
\frac{1}{2\sqrt{t}}\, \E[h'(\sqrt{1-t}F+\sqrt{t}Z)Z] = \frac{\sigma^2}{2}\, \E[h''(\sqrt{1-t}F+\sqrt{t}Z)].
\]
On the other hand, the relation $F= LL^{-1}F =-\delta DL^{-1}F$ and (\ref{ipp}) imply that
\begin{eqnarray*}
\frac{1}{2\sqrt{1-t}}\, \E[h'(\sqrt{1-t}F+\sqrt{t}S)F] &=& \frac{1}{2\sqrt{1-t}}\, \E[h'(\sqrt{1-t}F+\sqrt{t}Z)\delta(-DL^{-1} F)] \\
&=&  \frac{1}{2} \E[h''(\sqrt{1-t}F+\sqrt{t}Z)\langle DF, -DL^{-1} F\rangle_\HH].
\end{eqnarray*}
The conclusion follows from the fact that \[\big| \E[h(F)]- \E[h(Z)] \big|= \big| \Psi(1) - \Psi(0) \big|\leq \int_0^1 \big|\Psi'(t)\big|dt.\]

\noindent {\bf (2)} Here we follow the arguments contained in the proof of Theorem 3.1 in \cite{NP1}.
Define $h_\sigma(x) = h(\sigma x)$, $F_\sigma = \sigma^{-1}F$, and $Z_\sigma = \sigma^{-1}Z\sim N(0,1)$. Let $s$
be the solution of the Stein's equation associated with $h_\sigma$, i.e. $s$ solves the differential equation
$$ h_\sigma (x) - \E[h_\sigma(Z_\sigma)]=s'(x) - xs(x), \qquad x\in \R.$$
It is well-known that such a solution is given by $s(x) = \frac{1}{\varphi (x)} \int_{-\infty }^x (h_\sigma(t) - \Phi(t)) \varphi(t) dt$,
where $\varphi$ and $\Phi$ are the density and the distribution function of $N(0,1)$, respectively, and
$\| s'\|_\infty \leq \| h'_\sigma\|_\infty$. Since  $F_\sigma$ is a centered element of $\sk^{1,2}$ it holds that
$F_\sigma=LL^{-1}F_\sigma= -\delta D L^{-1}F_\sigma$. By integration by parts formula (\ref{ipp}) we deduce that
\begin{eqnarray*}
\big| \E[h(F)]- \E[h(Z)] \big| &=& \big| \E[h_\sigma(F_\sigma)]- \E[h_\sigma(Z_\sigma)] \big|\\
&=& \big| \E[s'(F_\sigma)- F_\sigma s(F_\sigma)] \big| \\
&=& \big| \E[s'(F_\sigma)(1 - \langle DF_\sigma,-DL^{-1}F_\sigma \rangle_\HH)] \big| \\
& \leq &\| s'\|_\infty \E| 1 - \langle DF_\sigma,-DL^{-1}F_\sigma \rangle_\HH| \\
& \leq &\| h'_\sigma\|_\infty \E| 1 - \langle DF_\sigma,-DL^{-1}F_\sigma \rangle_\HH|.
\end{eqnarray*}
We conclude by using the relations  $\| h'_\sigma\|_\infty = \sigma \| h'\|_\infty$ and \[\langle DF_\sigma,-DL^{-1}F_\sigma \rangle_\HH = \sigma^{-2}\langle DF,-DL^{-1}F\rangle_\HH.\]

\qed

When applied to the special case of a Gaussian random variable $F$, Theorem \ref{ProbApp} yields the following neat estimates.

\begin{cor}\label{CORProbApp}
Let $F\sim N(0,\gamma^2)$ and $Z\sim N(0,\sigma^2)$, $\gamma,\sigma>0$.
\begin{description}
\item[\bf (1)] For every $h$ twice continuously differentiable and with a bounded second derivative,
\begin{equation}\label{E : est11}
\big| \E[h(F)]- \E[h(Z)] \big| \leq \frac{\| h''\|_\infty}{2}\, |\sigma^2- \gamma^2 |.
\end{equation}
\item[\bf (2)] For every Lipschitz function $h$,
\begin{equation}\label{E : est22}
\big| \E[h(F)]- \E[h(Z)] \big| \leq \frac{\| h'\|_\infty}{\sigma\vee \gamma}\, |\sigma^2- \gamma^2 |.
\end{equation}
\end{description}
\end{cor}

\section{Proof of Theorem \ref{mainth}} \label{S : PROOFS}
\subsection{Preparation}
First, let us remark that
the process $X=(X_k)_{k\in \Z}$ can always be regarded as a subset of an isonormal Gaussian process
$\{W(u):~u\in \HH\}$, where $\HH$ is a separable Hilbert space with scalar product $\langle \cdot,\cdot \rangle_\HH$. More precisely, we shall assume (without loss of generality) that, for every $k\in \Z$ and every $1\leq l\leq d$, there exists $u_{k,l}\in\HH$ such that
\begin{equation*}
X_k^{(l)} = W(u_{k,l}), \qquad \mbox{and consequently} \qquad \langle u_{k,l}, u_{k',l'} \rangle_{\HH}=r^{(l,l')}(k-k'),
\end{equation*}
for every $k,k'\in \Z$ and every $1\leq l,l'\leq d$. Observe also that $\HH$ can be taken of the form $\HH = L^2(A,\mathscr{A},\mu)$, where $\mu$ is $\sigma$-finite and non-atomic.

Using the Hermite expansion (\ref{hermite}) of the function $f$ we obtain the Wiener chaos representation
\begin{equation*}
S_n = \sum_{m=q}^\infty I_m (g_m^n), \qquad g_m^n \in \HH^{\odot m},
\end{equation*}
where the kernels $g_m^n$ have the form
\begin{equation} \label{kernel}
g_m^n=\frac{1}{\sqrt{n}}
\sum_{k=1}^n \sum_{t\in \{1, \ldots, d\}^m} b_{t} \,u_{k,t_1} \otimes \cdots
\otimes u_{k,t_m}
\end{equation}
for certain coefficients $b_{t}$ such that the mapping $t\mapsto b_t$ is symmetric on $\{1,...,d\}^m$. One also has the identities
\begin{equation*}
\E[f_m^2(X_1)] = m! \sum_{t\in \{1, \ldots, d\}^m} b_{t}^2, \quad m\geq q, \qquad
\E[f^2(X_1)] = \sum_{m=q}^\infty m!\sum_{t\in \{1, \ldots, d\}^m} b_{t}^2.
\end{equation*}

Here is a useful preliminary result.

\begin{lem} \label{lem2} For the kernels $g_{m}^n$ defined in (\ref{kernel}) and any $1\leq e\leq m-1$ we obtain
the inequality, valid for every $n$,
\begin{equation} \label{kernelineq}
\|g_{m}^n \otimes_e g_{m}^n\|_{\HH^{\otimes 2(m-e)}}\leq \frac{d^m}{m!} \E[f_m^2(X_1)]\gamma_{n,m,e},
\end{equation}
where $\gamma_{n, m,e}$ is defined by (\ref{gamma}). Furthermore, we have that, for every $n$,
\begin{equation} \label{absch}
m!\|g_m^n\|_{\HH^{\otimes m}}^2 \leq \E[f_m^2(X_1)] (2K+d^q \theta),
\end{equation}
where the constants $K$ and $\theta$ are defined, respectively in (\ref{K}) and (\ref{theta}).
\end{lem}
\textit{Proof.} [Proof of (\ref{kernelineq})] Fix $1\leq e\leq m-1$. Observe that
\begin{eqnarray*}
g_{m}^n \otimes_e g_{m}^n &=&  \frac{1}{n}
\sum_{k_1,k_2=1}^n \sum_{t,s\in \{1, \ldots, d\}^m} b_{t} b_{s} \prod_{j=1}^e r^{(t_j,s_j)}(k_1 -k_2) \\[1.5 ex]
&\times& u_{k_1,t_{p+1}} \otimes \cdots
\otimes u_{k_1,t_m} \otimes u_{k_2,s_{p+1}} \otimes \cdots
\otimes u_{k_2,s_m}.
\end{eqnarray*}
We obtain
\begin{eqnarray*}
&&\|g_{m}^n \otimes_e g_{m}^n\|_{\HH^{\otimes 2(m-e)}}^2 \leq   \Big(\frac{d^m}{m!} \E[f_m^2(X_1)] \Big)^2
\\[1.5 ex]
&& \times n^{-2}\sum_{k_1, \ldots,k_4=1}^n \theta(k_1 - k_2)^e \theta( k_3 - k_4)^e \theta(k_1 - k_3)^{m-e} \theta( k_2 - k_4)^{m-e},
\end{eqnarray*}
where $\theta(j)$ is defined in (\ref{thetaJ}).
Since $\theta( k_3 - k_4)^e \theta(k_1 - k_3)^{m-e} \leq \theta( k_3 - k_4)^m + \theta(k_1 - k_3)^{m}$ we deduce
that
\begin{eqnarray*}
&&n^{-2}\sum_{k_1, \ldots,k_4=1}^n \theta(k_1 - k_2)^e \theta( k_3 - k_4)^e \theta(k_1 - k_3)^{m-e} \theta( k_2 - k_4)^{m-e}
\\[1.5 ex]
&&\leq 2n^{-1} \sum_{k\in \Z} \theta(k)^m \sum_{|k|\leq n} \theta(k)^e \sum_{|k|\leq n} \theta(k)^{m-e}\leq \gamma_{n,m,e}^2.
\end{eqnarray*}
Hence, we obtain (\ref{kernelineq}).

\noindent [Proof of (\ref{absch})] By the Cauchy-Schwarz inequality we have
$$m!\left|\Big \langle \sum_{t\in \{1, \ldots, d\}^m} b_{t} u_{k,t_1} \otimes \cdots
\otimes u_{k,t_m}, \sum_{t\in \{1, \ldots, d\}^m} b_{t} u_{k+l,t_1} \otimes \cdots
\otimes u_{k+l,t_m}   \Big \rangle_{\HH^{\otimes m}}\right| \leq \E[f_m^2(X_1)]$$
We  deduce, for any $m\geq q$,
\begin{eqnarray*}
&&m!\|g_m^n\|_{\HH^{\otimes m}}^2 = \frac{m!}{n}
\sum_{k_1,k_2=1}^n \sum_{t,s\in \{1, \ldots, d\}^m} b_{t} b_{s} \prod_{j=1}^m r^{(t_j,s_j)}(k_1 -k_2) \\[1.5 ex]
&&\leq \left( 2H \E[f_m^2(X_1)] + m! \sum_{|k|\geq K}  \theta(k)^m \Big(\sum_{t\in \{1, \ldots, d\}^m} |b_{t}| \Big)^2 \right)
\nonumber \\[1.5 ex]
&&\leq \E[f_m^2(X_1)] \Big(2K +  \sum_{|k|\geq K}  (d\theta(k))^m \Big)\leq \E[f_m^2(X_1)] (2K+d^q \theta) , \nonumber
\end{eqnarray*}
which implies (\ref{absch}). \qed \\ \\

The proofs of Point 1 and Point 2 in Theorem \ref{mainth} are similar, and are detailed in the subsequent two sections.

\subsection{Proof of Theorem \ref{mainth}-(1)}\label{SS : ProofP1}
First of all, we remark that $\theta (j)\rightarrow 0$ as $|j|\rightarrow \infty$, because $\sum_{j\in \Z} \theta (j)^q<\infty$.
This implies that $K<\infty$, where the constant $K$ is defined in (\ref{K}). Moreover, the asymptotic variance $\sigma^2$ is finite.
Indeed we have that $\sigma^2= \sum_{m=q}^\infty m!\|g_m^n\|_{\HH^{\otimes m}}^2 \leq \E[f^2(X_1)] (2K+d^q \theta)$ due to
(\ref{absch}). \newline \newline
The main proof is composed of several steps.

{\bf (a)} \textit{Reduction to a finite chaos expansion.} We start by approximating $S_n$ by a finite sum of multiple integrals.
Define
\begin{equation*}
S_{n,N} = \sum_{m=q}^N I_m (g_m^n).
\end{equation*}
Now, let
$h\in C^2 (\R)$ be a function with bounded second derivative. Since
\begin{equation*}
|h(x)-h(y)-h'(0)(x-y)|\leq \frac 12 \|h''\|_{\infty} (y-x)^2 + \|h''\|_{\infty} |x\|y-x|
\end{equation*}
for all $x,y\in \R$, we immediately obtain that
\begin{equation*}
|\E[h(S_n) ] -\E[ h(S_{n,N})]|\leq  \|h''\|_{\infty} \Big( \frac 12 \|S_n - S_{n,N}\|_{L^2(\PP)}^2
+  \|S_n \|_{L^2(\PP)} \|S_n - S_{n,N}\|_{L^2(\PP)}\Big).
\end{equation*}
By inequality (\ref{absch}) we deduce that
\begin{equation*}
\|S_n \|_{L^2(\PP)}^2 \leq (2K + d^q \theta ) \|f(X_1) \|_{L^2(\PP)}^2
\end{equation*}
and
\begin{eqnarray*}
&&\|S_n - S_{n,N}\|_{L^2(\PP)}^2 \leq (2K + d^q \theta )
\sum_{m=N+1}^\infty \|f_m(X_1) \|_{L^2(\PP)}^2  \\[1.5 ex]
&&\leq
(2K + d^q \theta ) \|f(X_1) \|_{L^2(\PP)}
\Big(\sum_{m=N+1}^\infty \|f_m(X_1) \|_{L^2(\PP)}^2 \Big)^{1/2},
\end{eqnarray*}
where $\theta$ is defined by (\ref{theta}). We conclude that
\begin{equation} \label{estimate1}
|\E[h(S_n) ] -\E[ h(S_{n,N})]|\leq \frac{3(2K + d^q \theta )}{2}  \|h''\|_{\infty} \|f(X_1) \|_{L^2(\PP)}
\Big(\sum_{m=N+1}^\infty \|f_m(X_1) \|_{L^2(\PP)}^2 \Big)^{1/2},
\end{equation}
which completes the first step. \qed
\\ \\
{\bf (b)} \textit{Bounds based on the interpolation inequality (\ref{E : est1}).} 
Let $Z_N$ be a centered Gaussian random variable with variance $\sum_{m=q}^N \sigma^2_m$. By using (\ref{E : est1}) in the special case $F=S_{n,N}, \, Z=Z_N$ and by applying e.g. (\ref{dtf}), we obtain
\begin{eqnarray} \notag
&& \E \Big[ h (Z_N ) ] -\E[ h(S_{n,N})\Big]\\
&& \leq \frac12 \|h''\|_\infty \Big\| \sum_{m=q}^N \sigma^2_m-\langle DS_{n,N},-DL^{-1}S_{n,N}\rangle_{\HH}
\Big \|_{L^2(\PP)}\notag \\
&&\leq \frac12 \|h''\|_\infty \sum_{p,s=q}^N \Big\|
\delta_{ps} \sigma_p^2-s^{-1}\langle DI_p(g_p^n),DI_s(g_s^n)\rangle_{\HH}
\Big \|_{L^2(\PP)},\label{estimate2}
\end{eqnarray}
where $\delta_{ps}$ is the Kronecker symbol. This completes the second step. \qed \\ \\
{\bf (c)} \textit{The final estimates.} Here we give the approximation of the term on the right-hand
side of (\ref{estimate2}). By (\ref{absch}) and the dominating convergence theorem we immediately
deduce that
\begin{equation*}
\E [s^{-1}\| DI_s(g_s^n)\|_{\HH}^2] =
s!\|g_s^n\|_{\HH^{\otimes s}}^2 \rightarrow  \sigma_s^2 = s!
\sum_{k\in \Z} \sum_{t,l\in \{1, \ldots, d\}^s} b_{t} b_{l} \prod_{j=1}^m r^{(t_j,l_j)}(k).
\end{equation*}
As in the proof of (\ref{absch})  we conclude that (recall that we assumed $n>K$)
\begin{eqnarray*}
|\E [s^{-1}\| DI_s(g_s^n)\|_{\HH}^2] - \sigma_s^2|
&\leq& s!  \sum_{|k|<K} \Big| \sum_{t,l\in \{1, \ldots, d\}^s} b_{t} b_{l} \prod_{j=1}^m r^{(t_j,l_j)}(k) \Big| \frac{|k|}{K} \\[1.5 ex]
&+& s!  \sum_{K\leq |k|<n} \Big| \sum_{t,l\in \{1, \ldots, d\}^s} b_{t} b_{l} \prod_{j=1}^m r^{(t_j,l_j)}(k) \Big| \frac{|k|}{K} \\[1.5 ex]
&+& s!  \sum_{|k|\geq n} \Big| \sum_{t,l\in \{1, \ldots, d\}^s} b_{t} b_{l} \prod_{j=1}^m r^{(t_j,l_j)}(k) \Big|  \\[1.5 ex]
&\leq & \E[f^2_s(X_1)] \left\{\frac{2K^2}{n}+d^q\left(
\sum_{|j|\leq n} \theta(j)^q \frac{|j|}{n} +\sum_{|j|>n}\theta(j)^q
\right)\right\}.
\end{eqnarray*}
Thus we have
\begin{eqnarray*}
&&\sum_{s=q}^N \Big\|
\sigma_s^2-s^{-1}\| DI_s(g_s^n)\|_{\HH}^2
\Big\|_{L^2(\PP)} \notag\\[1.5 ex]
&& \leq \sum_{s=q}^N \Big(|\E [s^{-1}\| DI_s(g_s^n)\|_{\HH}^2] - \sigma_s^2|
+ \sqrt{{\rm Var} [s^{-1}\| DI_s(g_s^n)\|_{\HH}^2]} \Big) \nonumber \\[1.5 ex]
&&\leq \E[f^2_s(X_1)] \left\{\frac{2H^2}{n}+d^q\left(
\sum_{|j|\leq n} \theta(j)^q \frac{|j|}{n} +\sum_{|j|>n}\theta(j)^q
\right)\right\} \nonumber \\[1.5 ex]
&&\qquad\qquad\qquad\qquad\qquad\qquad\qquad\qquad+ \sum_{s=q}^N \sqrt{{\rm Var} [s^{-1}\| DI_s(g_s^n)\|_{\HH}^2]} \nonumber \\[1.5 ex]
&&\leq 2(A_{1,n} + A_{3,n,N}),
\end{eqnarray*}
where the last inequality follows directly from Lemma \ref{lem1} and \ref{lem2}. On the other
hand, using the approximation (\ref{absch}) and Lemma \ref{lem1} and \ref{lem2} once again
we deduce that
\begin{eqnarray*}
&&\sum_{q\leq p<s\leq N} \Big\|
s^{-1}\langle DI_p(g_p^n),DI_s(g_s^n)\rangle_{\HH}
\Big\|_{L^2(\PP)} + \sum_{q\leq s<p\leq N} \Big\|
s^{-1}\langle DI_p(g_p^n),DI_s(g_s^n)\rangle_{\HH}
\Big\|_{L^2(\PP)} \nonumber \\[1.5 ex]
&&= \sum_{q\leq p<s\leq N} \frac{p+s}{p}~\Big\|
s^{-1}\langle DI_p(g_p^n),DI_s(g_s^n)\rangle_{\HH}
\Big\|_{L^2(\PP)}  \leq 2 (A_{4,n,N} + A_{5,n,N}).
\end{eqnarray*}
Thus, we conclude
\begin{equation} \label{estimate3}
\Big|\E \Big[ h \Big(\sum_{m=q}^N \sigma_m Y_m \Big) - h(S_{n,N})\Big]\Big|
\leq  \|h''\|_\infty (A_{1,n} + A_{3,n,N} + A_{4,n,N} + A_{5,n,N}),
\end{equation}
which finishes this step. \qed \\ \\
{\bf (d)} \textit{Putting things together.} We have the inequality
\begin{eqnarray}
|\E[h(S_n)] -\E[h(S)]|\label{UUU}
&\leq&
|\E[h(S_n)] -\E[h(S_{n,N})]|
+
\Big|\E[h(S_{n,N}) ] -\E[ h(Z_N)]\Big|\\[1.5 ex]
&&
+
\Big|\E[h(Z_N)] -\E[h(S)]\Big|.\notag
\end{eqnarray}
We are left with the derivation of a bound for the third term. By using (\ref{E : est11}) in the special case $F=Z_N$, we deduce that
\begin{equation}
\Big|\E[h(Z_N )] -\E[h(S)]\Big| \leq
\frac {1}{2} \|h''\|_\infty \sum_{m=N+1}^\infty  \sigma_m^2\leq \frac {2H+d^q \theta}{2}
 \|h''\|_\infty \sum_{m=N+1}^\infty  \E[f^2_{m} (X_1)],
\end{equation}
where the last inequality is deduced by Lemma \ref{lem2}. The latter is smaller than $\frac{1}{4} \|h''\|_\infty A_{2,N}$,
which together with (\ref{estimate1}), (\ref{estimate2}) and (\ref{estimate3}) completes the proof of Theorem
\ref{mainth}-(1). \qed
\subsection{Proof of Theorem \ref{mainth}-(2)}\label{SS : ProofP2}
Take $h$ Lipschitz and write the inequality (\ref{UUU}). Similar to (\ref{estimate1}) standard computations yield
\[
|\E[h(S_n)] -\E[h(S_{n,N})]|\leq \|h'\|_\infty \|S_n -S_{n,N}\|_{L^2(\PP)}\leq \|h'\|_\infty (2K + d^q \theta )^{1/2}\!\!
\sqrt{\sum_{m=N+1}^\infty \|f_m(X_1) \|_{L^2(\PP)}^2}.
\]
By applying inequality (\ref{E : est2}) in the case $F=S_{n,N},\, Z=Z_N$, one has that
\begin{eqnarray*}
&& \Big|\E[h(S_{n,N})] -\E[ h(Z_N)]\Big|\\
&& \leq \frac{\|h'\|_\infty}{\left(\sum_{m=q}^N\sigma^2_m\right)^{1/2}} \sum_{p,s=q}^N \Big\|
\delta_{ps} \sigma_p^2-s^{-1}\langle DI_p(g_p^n),DI_s(g_s^n)\rangle_{\HH}
\Big \|_{L^2(\PP)} \\
&&\leq \frac{\|h'\|_\infty}{\left(\sum_{m=q}^N\sigma^2_m\right)^{1/2}}2 (A_{1,n}+A_{3,n,N}+A_{4,n,N}+A_{5,n,N}),
\end{eqnarray*}
where the last inequality is obtained by reasoning as in Part (c) of Section \ref{SS : ProofP1}. Finally, an application of inequality (\ref{E : est22}) in the case $F=Z_N$ yields (since $\theta \geq 1$)
\[
\Big|\E[h(Z_{N}) ] -\E[ h(S)]\Big|\leq \frac{\|h'\|_\infty}{\sigma}\sum_{m=N+1}^\infty\sigma^2_m\leq  \frac{\|h'\|_\infty}{2\sigma}A_{2,N}.
\]
Putting the above estimates together yields the desired conclusion. \qed

\subsection{Proof of Theorem \ref{mainth}-(3)}\label{SS : ProofP3}
From
 \cite[Theorem 3.1]{Chen_Shao_sur},
 one can deduce that
 \[
 d_{Kol}(S_n,S)\leq 2\sqrt{d_W(S_n/\sigma,S/\sigma)}.
 \]
 Hence, we get the desired conclusion by combining this inequality with (\ref{mainestimate2}). \qed

\hskip2cm
 
\noindent
{\bf Acknowledgments}. We thank Arnaud Guillin for pointing us out an inaccuracy 
in a previous version.

\end{document}